%% file: taxcon.tex
\documentclass[11pt]{article}

\input{preamble}
\begin{document}


\input{title}


\section{Introduction}
\label{sec-intro}

This paper focuses on the feasible set $\Omega \subset \R^n$ of the 
general optimization problem
\begin{equation}
      \min_{x \in \Omega}  f(x),
     \label{eq-pb}
\end{equation}
where $f : \R^n \rightarrow \R \cup \{ \infty \}$ denotes an
extended-value objective function.
We propose a taxonomy of constraints, denoted \qrak,
whose development is motivated by the field of
{\em derivative-free optimization} (DFO),
and more precisely
  {\em black-box optimization} (BBO) and {\em simulation-based optimization} (SBO).
In BBO/SBO, the objective function $f$ and/or some constraints defining an instance 
of $\Omega$ are, or can depend on, the outputs of one or more {\em black-box simulations}.
We assume that SBO is the more general term; hence, we use it in 
the title of this work.
In typical settings, 
  evaluating the 
  simulation(s) is the primary bottleneck for an
  optimization algorithm; the time required to evaluate algebraic 
terms associated with other constraints or the objective is inconsequential 
relative to the time required to evaluate the
simulation 
components. 
In addition, 
simulations
may sometimes fail
to return a value, even for points inside $\Omega$.

Our taxonomy
 addresses a specific instance (or ``description'') 
of $\Omega$. This instance, rather than the mathematical problem 
(\ref{eq-pb}), will be passed to an optimization  solver (which may do 
some preprocessing of its own and then tackle a different instance).

To illustrate the distinction between problem and instance, we consider the 
two-dimensional linear problem
\begin{equation}
      \min_{x \in \R^2} \left\{x_1+x_2 : x_1\geq 0, ~x_2\geq 0\right\}.
          \label{eq:exprob}
\end{equation}
In fact, many instances of the feasible set $\Omega$ share 
a solution set with (\ref{eq:exprob}). For example, a different description can 
yield the same feasible set, either by chance, 
\[\Omega_1 = \{x\in \R^2: x_1\geq 0, ~x_1 x_2\geq 0\},\]
or as a result of some redundancy, 
\[\Omega_2 = \{x\in \R^2: x_1\geq 0, ~x_2\geq 0, ~2x_1+x_2\geq 0 \}.\]
Or, the feasible sets can differ from instance to instance, but 
the minimizers of $f$ over the sets are the same, whether indirectly,
\[\Omega_3 = \{x\in \R^2: x_1\geq 0, ~ x_1+2x_2\geq 0,~ x_1+x_2\leq 1\},\]
or explicitly,
\[\Omega_4 = \{x\in \R^2: x_1=0,~  x_2= 0\}.\]

In situations similar to these examples, one likely expects that a 
modern solver or modeling language---or 
even more classical techniques such as Fourier-Motzkin elimination---would 
perform preprocessing that would 
address redundancies, inefficiencies, and the like before invoking the heaviest 
machinery of a solver.
However, when the problem involves some black-box 
or simulation 
component, the situation, and hence such preprocessing, can be 
considerably more difficult.

More generally, the proposed classification is not
absolute: it depends on  
the entire set of constraint models specified in the instance
 and on the information
that the problem/simulation designer gives. 
For example, a simple bound constraint 
may be indicated as the output of a black box rather than expressed
algebraically,
leading to two different classes in the taxonomy.
Other examples of different constraints changing class
will be described after the taxonomy has been introduced.

Formally, we assume that a finite-dimensional instance $\Omega$ is specified by 
a collection of equations, inequalities, and sets:
\begin{equation}
 \Omega = \left\{ x\in \R^n : 
 c_i(x) = 0, \forall i\in \mathcal{I}; ~\;
 c_j(x) \leq 0, \forall j\in \mathcal{J}; ~\;
 c_k(x) \in \mathcal{A}_k, \forall k\in \mathcal{K} \right\},
 \label{eq-omega}
\end{equation}
where $\mathcal{I}, \mathcal{J}, \mathcal{K}$ are finite and
possibly empty. 
Semi-infinite problems can be 
treated by such a taxonomy but are not specifically addressed in this paper. 
Similarly, multi-objective optimization problems are easily encapsulated in our 
taxonomy but are not discussed specifically.
Note that the form of $\Omega$ in \eqref{eq-omega} is general enough to include
cases when a variable changes the total number of decision variables (such as
when determining the number of bus stations to build as well as their
locations).

As underscored in the recent book by Conn et al.~\cite{Dfobook}, derivative-free
optimization in the presence of general constraints has not yet been fully
addressed in the algorithmic literature or in benchmark papers such
as~\cite{MoWi2009} or~\cite{RiSa2010}. Even in broader SBO fields such as
simulation optimization and PDE-constrained optimization, a
disconnect often exists between what algorithm designers assume about a 
simulation and what
problem/simulation designers provide. 
In these communities, many different terms coexist for the same concepts, and
  unification is needed.
The proposed taxonomy of constraints
consolidates many previous terms
such as
{\em soft, virtual, hard, hidden, difficult, easy, 
open, closed}, and {\em implicit}.
Its purpose is to introduce a common language in order to
 facilitate dialog between algorithm developers,
optimization theoreticians, software users, and
application scientists formulating problems.

%

The paper is structured as follows. 
Section~\ref{sec-classes} presents
the taxonomy of constraints for SBO
and describes the different classes.
It also illustrates the taxonomy
through practical examples and situations.
Section~\ref{sec-lit-rev}
 puts the taxonomy in perspective with the existing literature.
 Putting the literature review toward the end of the paper here is deliberate 
and eases the presentation.
Section~\ref{sec-discussion} summarizes our contributions and discusses  
extensions to the taxonomy.

\section{Classes of constraints}
\label{sec-classes}

\begin{figure}[tb]
 \begin{center}
  \includegraphics [width=\textwidth]{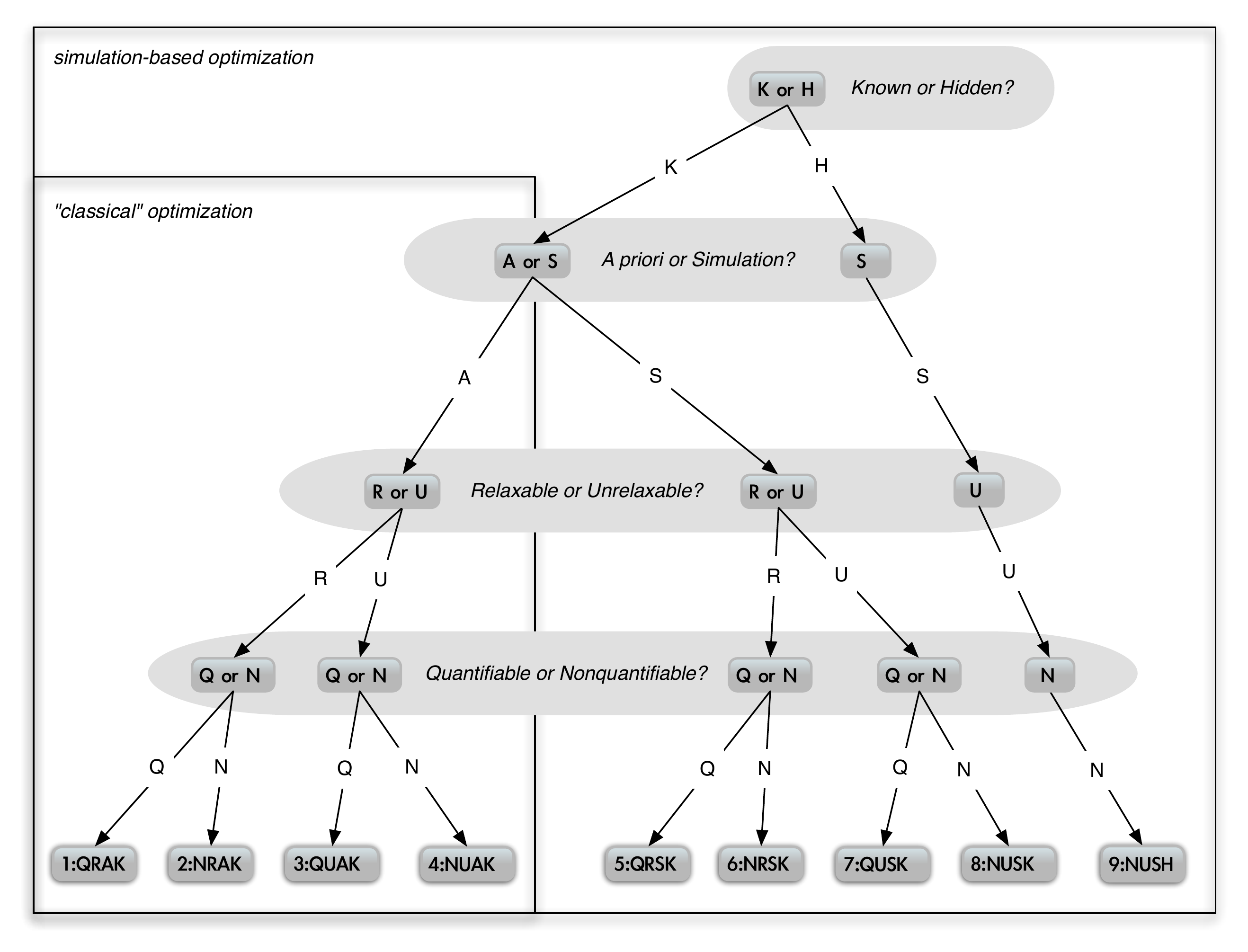}
   \end{center}
\caption{Tree-based view of the \qrak taxonomy of constraints.
Each leaf corresponds to a class of constraints.}
\label{fig-tree}
\end{figure}

\begin{figure}
 \begin{center}
  \includegraphics [width=12cm]{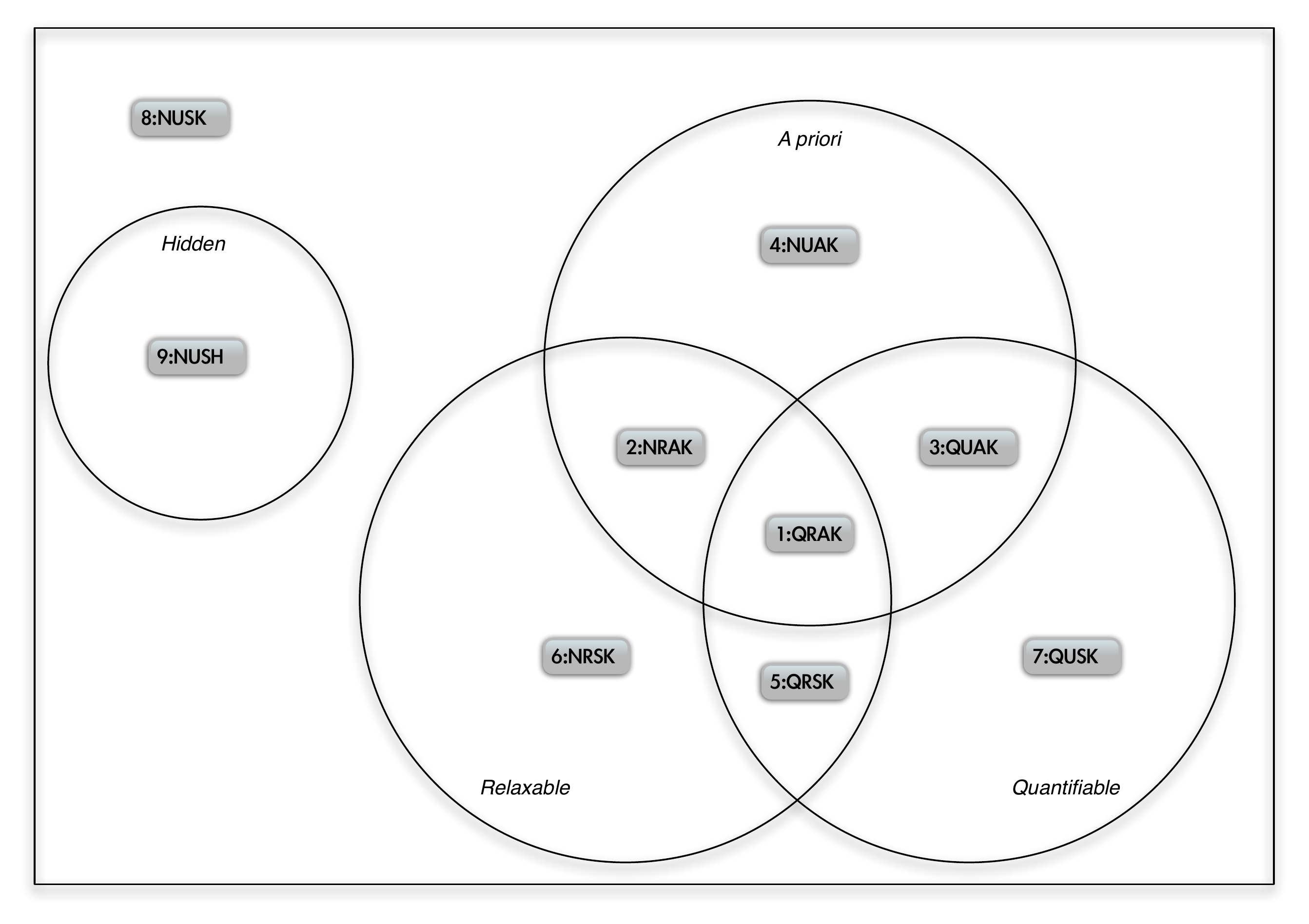}
 \end{center}
\caption{Venn diagram of the taxonomy of constraints. Each region corresponds to
a leaf in the tree of Figure~\ref{fig-tree}.}
\label{fig-Venn}
\end{figure}

This section introduces the \qrak taxonomy, which we present graphically by the 
tree of Figure~\ref{fig-tree}. 
An alternative and equivalent representation of the taxonomy using the same 
notations 
is given in the Venn diagram of Figure~\ref{fig-Venn}.

The letters defining the acronym of the taxonomy correspond to four types of left 
branches in the tree:
{\sf Q} is for {\sf Q}uantifiable,
{\sf R} is for {\sf R}elaxable,
{\sf A} is for {\sf A} priori,
and 
{\sf K} is for {\sf K}nown.
The corresponding right branches are identified with
{\sf N} for {\sf N}onquantifiable,
 {\sf U} for {\sf U}nrelaxable,
 {\sf S} for {\sf S}imulation,
 and
  {\sf H} for {\sf H}idden.

Each leaf of the tree in Figure~\ref{fig-tree} is identified with
a sequence of four letters, each entry taking one of two possible values.
The acronym of a leaf reads from the bottom to the root of the tree. 
As we argue later, not all 16 possible combinations of these letters are 
captured in the taxonomy, because hidden constraints take a special 
form. The nine possible constraint classes in the taxonomy are summarized in 
Table~\ref{fig-table}.

The two top levels of the tree are specific
  to SBO while the lower two are more general.
  In addition, most of constraints found in traditional nonlinear optimization 
(NLO) exist in the leftmost leaf.
In fact, general difficulty grows from left to right,
which outlines a preference for practitioners to model constraints
such that they appear in the most possible left part of the tree.
Further subdivisions (convexity, nonlinearity, etc.) are also important but
more focused on the NLO case and hence not discussed here.

Every constraint in an SBO problem instance fits in one leaf of the tree.
However, a constraint type from a classification scheme 
different from \qrak (e.g., bound constraint, nonlinear equality constraint) 
can  
correspond to several \qrak leaves at once.
In this case, we use the generic wildcard notation ``{\sf *}.''
For example, depending on the context, a bound constraint can be relaxable or unrelaxable.
It is clearly, however, a constraint that is known, a priori, and quantifiable.
In this case, the bound constraint is identified by {\sf Q*AK}.
%
The wildcard is not systemically used when the sense is obvious: For example,
we simply write {\sf S} instead of {\sf **S*}.
These issues will appear as natural as
we proceed with examples and formal definitions of each class/level of the tree, starting 
from the bottom and moving to the top.

\begin{table}
\center{
\caption{The taxonomy as a table where each column corresponds to a leaf in the
tree of Figure~\ref{fig-tree}
and to an intersection of regions in Figure~\ref{fig-Venn}. \label{fig-table}}
\begin{tabular}{|c|l|}
\hline
Leaf Number   & Name in  the \\
in Figure~\ref{fig-tree} & Taxonomy \\
\hline
\hline
1 & {\sf QRAK}  \\
\hline
2 & {\sf NRAK}  \\
\hline
3 & {\sf QUAK} \\
\hline
4 & {\sf NUAK} \\
\hline
5 & {\sf QRSK}  \\
\hline
6 & {\sf NRSK}  \\
\hline
7 & {\sf QUSK} \\
\hline
8 & {\sf NUSK} \\
\hline
9 & {\sf NUSH} (hidden) \\
\hline
\end{tabular}
}
\end{table}

\subsection[Quantifiable (Q) versus nonquantifiable (N)]
        {Quantifiable ({\sf Q})  versus nonquantifiable ({\sf N})}
For a nonquantifiable constraint, one has only a binary indicator saying
whether the constraint has been satisfied or violated. 
Consequently, an alternative term for such constraint is a {\em binary}
or {\em 0-1} constraint, but this does not have a natural
complementary term. Similarly, we avoid the terms {\em
measurable/nonmeasurable} in order to avoid confusion with {\em
measurable} in analysis. 

{\bf Definition:}
\label{def:qn}
A \underline{quantifiable constraint} is a constraint for which the degree of
feasibility and/or violation can be quantified. A \underline{nonquantifiable}
\underline{constraint} is one for which the degrees of satisfying or violating the
constraint are both unavailable.

The definition of a quantifiable constraint does not guarantee that
measures of both feasibility and violation are available. In particular, both of
the following are examples of quantifiable constraints.
\begin{description}
\item[Quantifiable feasibility:] {\em The time required for the underlying
simulation code to complete should be less than 10 seconds.}

Here, we have
access to the time that it took for the code to complete (and hence we know how
close we are to the 10-second limit), but the execution is interrupted if
it fails to complete within 10 seconds (and hence we will never know the degree to
which the constraint was violated). 
\item[Quantifiable violation:] {\em A time-stepping simulation should run to
completion (time $T$).} 

If the simulation stops at time $\hat{t}<T$, then
$T-\hat{t}$ measures how close one was to satisfying the constraint. 
\end{description}

A constraint for which both the degrees of feasibility and of violation are
available can be referred to as {\em fully quantifiable}. 

From a method or solver point of view, the distinction
between {\sf Q} and {\sf U} clearly is important. For example, if one wants to 
build
a model of the constraint,
{\sf Q} might imply {\em interpolation} whereas
{\sf U} might imply {\em classification}.

\subsection[Relaxable (R) versus unrelaxable (U)]
                  {Relaxable ({\sf R}) versus unrelaxable ({\sf U})}
\label{sec-relaxable}
The next notion addressed by the taxonomy is that of relaxability. 

{\bf Definition:}
 \label{def:ru}
A \underline{relaxable constraint} is a constraint that does not need
  to be satisfied
  in order to obtain {\em meaningful} outputs from the simulations
  in order to compute the objective and
  the constraints.
 An \underline{unrelaxable constraint}
   is one that must be satisfied for 
   meaningful outputs to be obtained.

In this definition, {\em meaningful} simulation output(s) means
that the values 
can be trusted as valid by an optimization algorithm 
and rightly interpreted when observed in a solution.

Typically, relaxable constraints are not part of a physical model 
  but instead represent some customer specifications 
  or some desired restrictions
  on the outputs of the simulation, such as a budget or 
  a weight limit.

Within an optimization method,
the implication regarding this {\sf R} versus {\sf U} property
is that
 all the iterates
  must satisfy unrelaxable constraints,
  while 
  relaxable constraints need be satisfied only
  at the proposed solution.
Said differently, infeasible points may be considered
as intermediate (approximate) solutions.

Alternative terms include
{\em soft} versus {\em hard}, {\em open} versus {\em closed}, and
{\em violable} versus {\em unviolable}; but these terms are often overloaded,
as we note in Section~\ref{sec-lit-rev}.

\subsection[A Priori (A) versus simulation-based (S)]
                  {A Priori ({\sf A}) versus simulation-based ({\sf S})}
A simulation constraint is specific to BBO/SBO.
The nature of a simulation constraint is such that a potentially costly
call to a computer simulation must be launched in
order to evaluate the constraint. We
note, however, that this constraint evaluation may not ultimately prove to be
costly. For example, the simulation could include
a constraint that is cheap to evaluate and can be used
as a flag to avoid any further computation;
such a constraint is still defined by our taxonomy to be an {\sf S}
constraint (more specifically, a {\sf *USK} constraint).

{\bf Definition:}
 \label{def:as}
 An \underline{a priori constraint} is a constraint
 for which feasibility can be confirmed without running
 a simulation.
A \underline{simulation-based constraint}
(or simulation constraint)
  requires running a simulation
  to verify feasibility.

Simple examples of a priori constraints include one-sided bounds
and linear equalities. 
A Priori constraints, however, can include very general and special
formulations,
such as
semidefinite programming constraints,
constraint programming constraints (e.g., {\em all different}, {\em
ordered}),
or some constraints relative to the nature of the variables: reals or integers or binary or
categorical.

One can easily appreciate that a solver should want to evaluate {\sf *UA*}
(unrelaxable, a priori) constraints first
and avoid a simulation execution
if the candidate is infeasible---especially when the simulation
is costly.
For {\sf *RA*} (relaxable, a priori) constraints, it is not as clear whether an 
algorithm would benefit from a similar ordering of constraint evaluations. 
For example, should noninteger input values be passed to a simulator that
may then end up rounding to the nearest integer within the simulation?
The answer depends on the context.

An alternative to ``simulation" is {\em a posteriori}~\cite{AuDaOr13a};
alternative terms for ``a priori" include
{\em algebraic} or {\em algebraically available},
{\em analytic},
{\em closed-form}, 
{\em expressible}, and
{\em input-constraint}.
An {\em algebraic function} is
usually defined to be one that satisfies an equation that can be expressed as a
finite-degree polynomial with rational coefficients. Unfortunately this 
definition does
not include transcendental functions ($e^x$, etc.). 
Some modeling languages, such as \gams~\cite{BrKeMe88a}, 
already use this term (\gams is ``generalized algebraic" to include
available transcendentals).
Formally, an {\em analytic function} is usually one that locally has a
convergent power series; this rules out simple nonsmooth functions. 
The idea behind the term
{\em input-constraint} is that {\sf A} constraints can be seen
as simply related to the 
inputs $x$, whereas {\sf S} constraints are somehow
   expressed as a function of the simulator.

\subsection[Known (K) versus hidden (H)]
                  {Known ({\sf K}) versus hidden ({\sf H})}
\label{sec-hidden}
The final distinction in the taxonomy is specific to BBO/SBO.

{\bf Definition:}
 \label{def:ku}
A \underline{known constraint}
is a constraint that is explicitly given in the problem formulation.
A \underline{hidden constraint} is not explicitly known to the solver.

The majority of constraints that one encounters when solving SBO
problems --- especially when an optimizer is involved early in the modeling and
problem formulation process --- are {\em known} to the optimizer. 
A hidden constraint typically (but not necessarily) appears when
the simulations crashes.
For such constraints, we can detect only violations,
typically when some error flag or exception is raised.
However, a violation may go unnoticed.
Alternative terms include
{\em Unknown}, {\em Unspecified}, and {\em Forgotten}.

A hidden constraint is not necessarily a bug in the simulator.
For example, consider the problem $\min \{f(\log x) : x \in \R \}$ with $f$ some
simulation-based function from $\R$ to $\R$.
If the constraint $x > 0$ is expressed in the description of the problem, then
it is an a priori constraint. Otherwise, it is hidden and 
 can be observed only for negative or null values of $x$.
This constraint may have been stated explicitly inside the simulator in order to
avoid a crash and raise some
flag; but as far as it not indicated to the solver, it remains {\sf H}.

As shown in Figure~\ref{fig-tree}, the {\sf H} branch of the tree is the
only one that goes directly to a terminal leaf.
A hidden constraint cannot be a priori (by definition)
and quantifiable 
(we do not know what to quantify).
A hidden constraint  also cannot be relaxable since
the violation/satisfaction cannot be detected if the outputs are always meaningful.

Note that the boundary between a hidden ({\sf NUSH}) constraint and a {\sf
NUSK} constraint is thin.
In the {\sf NUSK} case, however, the constraint is explicitly given, and its satisfaction
can be checked.
These subtle differences are emphasized in the presence of several different
hidden constraints:
When the simulation crashes, one has no way of knowing exactly what went wrong,
a situation that would have been different
if these constraints had been expressed with flags by the modeler.

\subsection{Short case studies}
\label{sec-examples}

The previous examples were related to the four levels of
decision in the taxonomy.
We now show that each of the nine leaves of the tree in Figure~\ref{fig-tree}
(similarly,
each row of Table~\ref{fig-table}) is nonempty, and we illustrate some
situations that belong to each leaf.

\begin{enumerate}
 
\item[1:] 
{\sf QRAK} (Quantifiable Relaxable A Priori Known):
Probably the most common type of constraint found in classical nonlinear optimization.

\begin{itemize}

\item
$\sum\limits_{i=1}^n x_i \leq 100$: If each $x_i$ represents an amount of money,
 this constraint defines a budget.

\item
Relaxable discrete variable:
$x_i\in \{0,1\}$ for some indices $i$.
Then $\min\bigl\{|x_i|,$  $|1-x_i| \bigr\}$ provides the violation measure.

\end{itemize}

 \item[2:]
{\sf NRAK} (Nonquantifiable Relaxable A Priori Known):
A good example is a categorical variable constrained to a subset of its  
possible values:

\begin{itemize}

\item
A simulator can work in two modes depending
  on the value of a binary flag $x \in \{0,1\}$.
 If $x=1$, the simulation is costly but more accurate.
  If $x=0$, it is cheap but imprecise.
  We want a solution that has been validated with $x=1$, but an optimization
  algorithm can set $x=0$ at intermediate points.
  The {\sf NRAK} constraint is ``$x=1$".

\item
Consider a simulator that drives a {\sf C++} compilation,
with the two categorical variables
$x_1 \in \{ \text{\sf gcc}, \text{\sf icc} \}$
and $x_2 \in \{ \text{\sf O2} , \text{\sf O3} \}$.
We want the final solution to have $x_2=\text{\sf O2}$
  if $x_1=\text{\sf gcc}$, and these two constraints are of type {\sf
NRAK}. (Note that the two set constraints, $x_1 \in \{ \text{\sf gcc}, \text{\sf
icc} \}$ and $x_2 \in \{ \text{\sf O2} , \text{\sf O3} \}$, may be {\sf NUAK};
see below.)


\end{itemize}

 \item[3:]
{\sf QUAK} (Quantifiable Unrelaxable A Priori Known):

\begin{itemize}


\item
Well rates in groundwater problems:
If we can simulate only extraction (and not injection), then the
constraints $r_i \geq 0$ for all well indexes $i$ are of type {\sf QUAK}.

\item
Decision variables must be ordered or be all different.

\end{itemize}

 \item[4:]
{\sf NUAK} (Nonquantifiable Unrelaxable A Priori Known):
Categorical variables are a typical example:
$ \text{\sf compiler}  \in \{ \text{\sf gcc} , \text{\sf icc} \}$.

 \item[5:]
{\sf QRSK} (Quantifiable Relaxable Simulation Known):
Simply consider a requirement on one of the simulation's
output, such as the following:

\begin{itemize}

\item
A budget based on economical criteria, $S(x)\leq b$.


\item
In the context of optimization of algorithm parameters, 
the percentage of problems solved by the algorithm under consideration 
must be $100\%$.

\end{itemize}

 \item[6:]
{\sf NRSK} (Nonquantifiable Relaxable Simulation Known):

\begin{itemize}

\item A simulator displays a flag indicating whether a toxicity level
   has been reached during the simulation,
   but we know neither when this occurred nor the level of toxicity.
   
\item
  The simulator indicates whether the power consumption
  remained under 100W, but we have access only to the notification.

\end{itemize}

 \item[7:]
{\sf QUSK} (Quantifiable Unrelaxable Simulation Known):
One of the outputs $c_S(x)$ of the simulation
is a concentration level; if it is below zero, the simulation stops and
displays {\sf NaN} for all the outputs
except $c_S$.

 \item[8:]
{\sf NUSK} (Nonquantifiable Unrelaxable Simulation Known):

\begin{itemize}

\item
A flag indicates that the convergence of some specific and identified
numerical method inside the simulation could not converge.

\item
An error number/code with associated documentation is obtained.
\end{itemize}

These are not hidden constraints since the reason for the violation can be 
identified.
However, a single binary flag indicating that the simulation failed is considered as a hidden constraint.
In the same way, an error message that cannot be interpreted is equivalent to such a flag
and hence should be interpreted as hidden.

 \item[9:]
{\sf NUSH} (Hidden):

The simulation failed to complete and nothing is displayed,
 or a simple flag is raised or an undocumented error number indicated.

\end{enumerate}

\section{Literature review}
\label{sec-lit-rev}

In this section, we 
review
the existing literature and collect terminology from the BBO, DFO, and SBO 
communities
in order to 
unify and relate our taxonomy to past terms and formulations 
and to highlight inconsistencies among previous conventions.
This context also underpins the naming
conventions used in \qrak\ and the more formal definitions
on which the taxonomy is built.
Some of the terms from the literature
may have been used to define
an alternative classification of the constraints,
and some of them
have already been mentioned in our presentation,
such as {\em soft} versus {\em hard} in Section~\ref{sec-relaxable}.
We also
survey early uses of various terms (e.g., hidden
constraints) for a historical perspective,
and we illustrate the use of the taxonomy in the context of modeling languages, 
algorithms, and some
specific applications.

Before proceeding, we note that the proposed classification is not
related to the constraint programming field~\cite{RossiVB06},
where, within a specific context, constraints can be expressed as logical 
prepositions
treated by specialized algorithms.

\subsection{Hidden constraints}

The term {\em hidden constraint}  corresponds
to the {\sf NUSH} leaf in the tree of Figure~\ref{fig-tree}.
It often appears in the literature on derivative-free optimization.
In the modern literature, this term is typically attributed to Choi and 
Kelley~\cite{ChKe00a}, who say that 
a hidden constraint is ``the requirement that the objective be defined." 
This
definition is used in Kelley's implicit filtering software~\cite{Kelley2011} and
has been used to solve several examples (see, e.g.,~\cite{Carter2001,ChoiOPTE})
whereby a hidden constraint is said to be violated whenever flow conditions are
found that prevent a simulation solution from existing.
The term is also used by the authors of the {\sf SNOBFIT} 
package~\cite{Huyer2008} to capture 
when ``a requested function value may turn out not to be obtainable."
To handle such constraints, {\sf SNOBFIT} assigns an 
artificial value, based on the values of nearby points, to the points where 
such a constraint was violated.
A more recent reference to the term is in~\cite{GrVi2014},
where hidden constraints are
  ``constraints which are not part of the problem specification/formulation 
and their manifestation comes in the form of some indication that the objective 
function could not be evaluated.''

In fact, the term had been previously used in the context of optimization. The
earliest published instance of hidden constraint that we are aware of is from
1967~\cite{Avriel67} and involved optimizing the design of a condenser. In
this case, after a design was numerically evaluated, one needed to verify that
the Reynolds number obtained was high enough to justify use of the equations in
the calculations. 

\subsection{{\em Introduction to Derivative-Free Optimization} textbook}

Although the DFO book~\cite{Dfobook} focuses on unconstrained optimization, it 
uses
definitions of relaxable and hidden constraints similar to those  
used in \qrak.
The book mentions that
unrelaxable constraints
``have to be satisfied at all iterations" of an algorithm
while ``relaxable constraints need only be satisfied approximately or 
asymptotically.'' 
Moreover,
constraints for which derivatives are not available and which are typically 
given by a black box,
are denoted as
{\em derivative-free constraints}.
Although in general these constraints can be treated as relaxable,
some situations require them to be unrelaxable.
Doing so demands a feasible starting point, which may be difficult to 
obtain in practice.
Moreover,
hidden constraints are seen as an extreme case of such unrelaxable constraints.
They are defined as constraints that 
\begin{quote}
``are not part of the problem specification/formulation, and their manifestation 
comes in the form of
some indication that the objective function could not be evaluated.'' 
\end{quote}
The authors of~\cite{Dfobook} state that hidden constraints have historically 
been treated only by heuristic approaches or by the {\em extreme-barrier} 
approach, which uses extended-value functions in an attempt to establish 
feasibility.

Scheinberg et al.~\cite{CST98} refer to {\em virtual constraints} as
``constraints that cannot explicitly be measured." Only the
satisfiability 
of
such constraints can be checked, and this is assumed to be a computationally
expensive procedure. In our taxonomy such constraints are {\rd {\sf N**K}}.
%

\subsection{Unrelaxable and relaxable constraints}

The terms {\em unrelaxable} and {\em relaxable} are widely used in 
the literature
(see~\cite{Dfobook}). The related notions of {\em hard} 
and {\em soft} constraints appear almost as frequently but with several 
different meanings. Here, we follow the convention
of~\cite{griva2009linear}:
\begin{quote}
``To resolve this, the requirements are usually broken up into ``hard"
constraints for which any violation is prohibited, and ``soft" constraints for
which violations are allowed. Typically hard constraints are included in the
formulation as explicit constraints, whereas soft constraints are incorporated
into the objective function via some penalty that is imposed for their
violation."
\end{quote}
That is, we view soft constraints as being handled by either additional
objectives or additional objective terms.  A nice pre-1969 history of ways to
move constraints into the objective can be found in~\cite{fiacco1968nonlinear}.
Another example comes from {\sf SNOBFIT}~\cite{Huyer2008}, 
where soft constraints are
``constraints which need not be satisfied accurately."  
Other uses of hard and soft constraint can be found, for example,
 in~\cite{Griffin01012010}. There, the authors refer to soft constraints as those
``that need not be satisfied at every iteration," a definition that 
 is directly related to our term {\em unrelaxable}. 
A similar notion is used in~\cite{GrVi2014}:
``Relaxable constraints need only be satisfied approximately or 
asymptotically.'' But our definition requires that a solution satisfy relaxable 
constraints, and hence the degree of ``approximate'' satisfaction must be 
specified in the problem instance.
In~\cite{Martelli2014108}, constraints are divided into relaxable and 
unrelaxable constraints, where 
unrelaxable constraints
\begin{quote}
``cannot be violated by any considered solution because they guarantee 
either the successful evaluation of the black-box function \ldots or the 
physical/structural feasibility of the solution"
\end{quote}
and relaxable constraints 
``may instead be violated as the objective 
function evaluation is still successful."

\subsection{Modeling languages and applications}

Several modeling languages and collections of test problems,
such as
\ampl \cite{FoGaKe03a},
\cutest~\cite{cutest},
\gams~\cite{BrKeMe88a}, or
{\sf ZIMPL}~\cite{Koch2004},
use the following classic ways of categorizing constraints:
fixed variables; bounds on the variables;
adjacency matrix of a (linear) network;
linear, quadratic,
equilibrium, and
conic constraints;
logical constraints found in constraint programming; and 
equalities or inequalities.
Usually, the remaining constraints are qualified as ``general,''
a term frequently used in classical nonlinear optimization.
All these constraints fit as {\sf **AK} constraints in the 
``classical optimization" portion of the tree of Figure~\ref{fig-tree}.

The \qrak\ taxonomy can be illustrated on the following examples
of SBO problems from the recent literature.
The community groundwater problem~\cite{FoRe_etal08} has only bound constraints 
({\sf **AK}), while the LOCKWOOD problem~\cite{Matott2011} has a linear 
objective and simulation constraints ({\sf S}); 
different simulation-based instances of the LOCKWOOD constraints are 
considered in~\cite{AKSW12} alongside solution methodologies for the resulting 
formulations. 
The  STYRENE problem from~\cite{AuBeLe08} has 11 simulation constraints corresponding to Leaves~5
and~8 of the tree of Figure~\ref{fig-tree}: 
7 quantifiable and relaxable constraints {\sf QRSK},
and 4 unrelaxable binary constraints {\sf NUSK}.

\subsection{Algorithms and software for constrained problems}
\label{sec-software}

To motivate the opportunities that such a constraint taxonomy affords, 
we briefly describe how some algorithms 
and software address different types of constraints,
  using the taxonomy syntax.
  
In general, most general-purpose software packages consider {\sf QR*K} 
constraints, but some tend
to use exclusively algebraic forms (e.g., box, linear, quadratic, convex).
Furthermore, relaxable constraints often are also assumed to be  
quantifiable.
Several packages allow for a priori constraints, but some assume that these 
cannot be relaxed, while
others assume that they can.
 
The package {\sf SNOBFIT}~\cite{Huyer2008}
treats {\em soft} and also {\sf NUSH} (hidden) constraints. The software 
{\sf SID-PSM}~\cite{CuVi05web}
handles constraints with derivatives and {\sf U} (unrelaxable) constraints.
The {\sf DFO} code~\cite{DFO} (which we distinguish from the general class of 
optimization problems without derivatives) considers
{\sf NUSH} (hidden), {\sf NU*K}, and {\sf Q*AK} constraints.
On the {\sf DFO} solver page~\cite{DFO},
the authors
 recommend moving {\sf S} (simulation, {\em difficult}) constraints to the 
objective function, while keeping 
{\em easy} constraints (with derivatives) inside the trust-region subproblem; 
the authors also describe {\em virtual} constraints as {\sf N} 
(nonquantifiable) constraints and
  recommend using an extreme-barrier approach.
The {\sf HOPSPACK} package~\cite{Hops20-Sandia} explicitly addresses
integers; linear equalities and inequalities; and general inequalities and 
equalities.
Depending on the type of constraint, {\sf HOPSPACK} assumes that the
constraint is relaxable (e.g., general equality constraints) or unrelaxable
(e.g., integer sets). 
In \nomad~\cite{Le09b}, the progressive-barrier technique~\cite{AuDe09a}
is used for the {\sf QRSK} constraints,
and special treatment (such as projection) is applied for some {\sf Q*AK} 
constraints (i.e., bounds and integers).
The extreme barrier is used for all other constraints, including hidden 
constraints.

In PDE-constrained optimization, solution approaches can be loosely classified
into ``Nested Analysis and Design'' (NAND) and ``Simultaneous Analysis and
Design'' (SAND) approaches~\cite{BGHvBWbook}. In NAND approaches, the state
variables of the PDE constraints are not treated as decision (optimization)
variables and hence the solution of the PDE (for the state variables) is a
simulation constraint. This situation exists even if the simulation is not just
a black box, but also returns additional information (e.g,. sensitivities,
adjoints, tolerances). In the NAND approach, the state variables are included as
decision variables and hence the PDE reduces to a set of algebraic equations
(and therefore {\sf *A*K} constraints in our taxonomy).

\subsection{Other related work}

Previous classifications also have been proposed. An example is the 
mixed-integer 
programming 
 classification~\cite{Nemhauser1991} for linear inequalities, linear equations,  
continuous parameters, 
and discrete parameters.

The closest related work toward a more compete characterization of constraints 
is that of Alexandrov and Lewis~\cite{Alexandrov1999},
who examined different formulations for general problems arising in
multidisciplinary optimization (MDO). These authors considered constraint sets
partitioned along three axes: open (closed) disciplinary analysis, open (closed)
design constraints, and open (closed) interdisciplinary consistency
constraints. 
They showed that of the eight possible combinations, only four
were possible in practice. They referred to {\em closed} constraints as
those 
\begin{quote}
``assumed to be satisfied at every iteration of the optimization. 
If
the formulation does not necessarily assume that a set
of constraints is satisfied, we will say that that formulation is open with
respect to the set of constraints."
\end{quote}
This convention has subsequently been used
by others in the MDO community (see, e.g.,~\cite{Tosserams2009}).

The notion of unknown constraints appears in~\cite{GrLe2011b} 
but is not equivalent to its use in our taxonomy; rather, it corresponds 
to constraints given by a black box.
Note that the same authors and others addressed hidden constraints 
in~\cite{LeGrLiGr2011}.

Additional terms for describing general constraints are found in the 
literature.
For example,
{\em chance constraints}~\cite{ChCoSy1958} are
constraints whose satisfaction requirement depends on a probability.
{\em Side constraint} is
a generic term sometimes used to qualify constraints that are not lower or 
upper bounds 
or to distinguish new constraints added to a preexisting model; 
see~\cite{Aggarwal1982287} for an example.
Other terms include
the notions of {\em vanishing} constraints~\cite{AcKa2008},
{\em complementarity} constraints,
or {\em variational inequalities}~\cite{LuPaRa1996}.
Conn et al.\ describe {\em easy} constraints and {\em difficult} constraints 
as follows~\cite{CST98}:
\begin{quote}
``Easy constraints are the constraints whose values and derivatives can be 
easily computed,"
\\and\\
``Difficult constraints are constraints whose derivatives are not available 
and
whose values are at least as expensive to compute as that of the objective 
function.''
\end{quote}
The latter definition is similar to the {\em derivative-free} constraints 
described 
in~\cite{Dfobook}. This characterization differs from our 
proposed taxonomy, which does not seek to guarantee an ordering with 
regard to the computational expense of evaluating a constraint and/or 
establishing feasibility with respect to the constraint.

\section{Discussion}
\label{sec-discussion}

This work proposes a
unification of past conventions and terms into a single taxonomy, denoted \qrak, 
which targets the constraints encountered in
 simulation-based optimization.
The taxonomy  has an intuitive representation as a tree where each leaf 
describes
one of nine types of possible constraints.
In addition, examples have been given for each constraint type and their 
possible treatment in applications and algorithms. 

We propose that 
BBO, DFO, and SBO software and algorithms should adopt this taxonomy 
for two important reasons.
The first is unification, so that
researchers in the field use the same terms and
 practitioners and algorithm developers
share the same language.
The second reason is that the taxonomy is a tool
to better identify constraint types and thereby achieve effective 
 algorithmic treatment of more general types of constrained optimization 
problems.

Future work is related to extensions to the taxonomy.
One can refine the tree in Figure~\ref{fig-tree}, depending on the context,
by adding subcases to the leaves.
Such extensions within {\sf QRAK} could include
stochastic, convex, linear, and smooth constraints (i.e., constraints for which 
derivatives are available).
Equality, inequality, or set membership is also an option:
For example, an equality {\sf N*SK} constraint 
is difficult (impossible?) to treat,
whereas 
an equality {\sf Q*AK} constraint 
may be easy.
At a different level, we consider the addition
of three branches from each {\sf Q} node: quantifiable
feasibility only, quantifiable violation only, and fully quantifiable.
There is also a limit to being unrelaxable:
So far we say that a constraint is unrelaxable if it is unrelaxable
at some point, and we may want to specify such limits when they are known.
%

\section*{Acknowledgments}

This material was based upon work supported by the U.S.\ Department of Energy, 
 Office of Science, Office of Advanced 
Scientific Computing
Research,  under Contract
DE-AC02-06CH11357;
by the Air Force Office of Scientific Research Grant FA9550-12-1-0198;
and
by the Natural Sciences and Engineering Research Council of Canada Discovery Grant 418250.
The authors thank the American Institute of Mathematics
for the SQuaRE workshops that allowed them to initiate their discussions on the 
taxonomy; the authors thank the participants:
Bobby~Gramacy, Genetha~Gray, Herbie~Lee, and Garth~Wells.
Thanks also to Charles~Audet 
for his useful comments and suggestions.

\bibliographystyle{siam}
\bibliography{cited}

\end{document}

%% file: preamble.tex
\usepackage{hyperref}
\hypersetup{
    unicode = false,
    pdftoolbar = true,
    pdfmenubar = true,
    pdffitwindow = true,
    pdftitle = {Taxonomy of constraints},
    pdfauthor = {S. Le Digabel and S.M. Wild},
    pdfsubject = {Taxonomy},
    pdfnewwindow = true,
    pdfkeywords = {Taxonomy},
    colorlinks = true,
    linkcolor = blue,
    citecolor = blue,
    filecolor = black,
    urlcolor = blue,
    breaklinks = true
}
\usepackage{amsmath}            
\usepackage{amssymb}            
\usepackage{float}              
\usepackage{url}           	
\usepackage{xspace}		

\usepackage{multirow}
\newcommand{\R}{\mbox{${\mathbb R}$}}           

\newcommand{\qrak}{{\sf QRAK}\xspace}

\newcommand{\gams}{{\sf GAMS}\xspace}
\newcommand{\ampl}{{\sf AMPL}\xspace}
\newcommand{\nomad}{{\sf NOMAD}\xspace}

\newcommand{\cutest}{{\sf CUTEst}\xspace}

\newcommand{\SDNote}[1]           
{\textcolor{red}{#1}\marginpar[\textbf{$\Longrightarrow$}]
{\textbf{$\Longleftarrow$}}}

\newcommand{\SWNote}[1]           
{\textcolor{blue}{#1}\marginpar[\textbf{$\Longrightarrow$}]
{\textbf{$\Longleftarrow$}}}

\usepackage{epsfig}
\usepackage{pgf}

\definecolor{Red}{rgb}{1,0,0}
\definecolor{Gray}{rgb}{0.2,0.2,0.2}

\newcommand{\rd}{\color{black}}

%% file: title.tex

\clearpage

\begin {center}
{\Large\textbf{A Taxonomy of Constraints in Simulation-Based 
Optimization}}
\bigskip

\end {center}

\begin{center}
\textbf{S\'ebastien Le Digabel\footnotemark [1] 
and 
Stefan M. Wild\footnotemark[2]}
\end {center}

\footnotetext [1] { GERAD and D\'epartment de Math\'ematique et de G\'enie Industriel, \'Ecole Polytechnique de Montr\'eal, Montr\'eal, QC H3C 3A7, Canada.
\url{https://www.gerad.ca/Sebastien.Le.Digabel}.}
\footnotetext [2] { Mathematics and Computer Science Division, Argonne
  National Laboratory, Argonne, IL 60439, USA.
\url{http://www.mcs.anl.gov/~wild/}.
  }

\bigskip

{\bf Abstract:}
The types of constraints encountered in black-box and simulation-based
optimization problems differ significantly from those treated
in nonlinear programming.  We
introduce a characterization of constraints to address this situation.
We provide formal definitions for several constraint classes and present illustrative examples in the context of the
resulting taxonomy. This taxonomy, denoted \qrak,  is useful for modeling and problem
formulation, as well as optimization software development and deployment. 
It can also be used as the basis for a dialog with practitioners in moving problems to increasingly solvable branches of 
optimization.

{\bf Keywords:}
Taxonomy of constraints, Black-box optimization, Simulation-based optimization.
